\documentclass[12pt,twoside]{article}
\usepackage{amsmath,amsbsy,amsfonts}
\usepackage[]{graphicx}
\def\dis{\displaystyle}

\newtheorem{lem}{Lemma}[section]

\newtheorem{theo}{Theorem}[section]
\newtheorem{coro}{Corollary}[section]

\newtheorem{theorem}{Theorem}[section]
\newtheorem{lemma}[theorem]{Lemma}

\newcommand{\fim}{\hfill\rule{2mm}{2mm}}

\newcommand{\RR}{\mathrm{I\!R\!}}

\newcommand{\ds}{\displaystyle}
\def\dis{\displaystyle}

\let\Section=\section
\def\section{\setcounter{equation}{0}\Section}

\title{Multiple solutions for a NLS equation with critical growth and magnetic field.}

\author{Claudianor O. Alves\footnote{Supported by INCT-MAT, PROCAD,NPq/Brazil
620150/2008-4 and 303080/2009-4} \\
\noindent Universidade Federal de Campina Grande \\
\noindent Departamento de Matem\'atica e Estat\'istica\\
\noindent CEP: 58429-900, Campina Grande - Pb, Brazil.\\
\noindent e-mail: coalves@dme.ufcg.edu.br\vspace{0.5cm}  \\
\noindent Giovany M. Figueiredo \thanks{Supported by CNPq/PQ 300705/2008-5 }\\
\noindent Universidade Federal do Par\'a \\
\noindent Faculdade de Matem\'atica \\
\noindent CEP: 66075-110 Bel\'em - Pa , Brazil. \\
\noindent e-mail: giovany@ufpa.br \vspace{0.5cm}\\
}
\vspace{0.5cm}

\date{}
\usepackage{color}

\begin{document}

\maketitle

\begin{abstract}
In this paper, we are concerned with the multiplicity of nontrivial
solutions for the following class of complex problems
$$
(-i\nabla - A(\mu x))^{2}u= \mu |u|^{q-2}u + |u|^{2^{*}-2}u \ \mbox{in}
\ \Omega , \ \ \ \
          u \in H^{1}_{0}(\Omega,\mathbb{C}),
$$
where $\Omega \subset \mathbb{R}^{N} ( N \geq 4 )$ is a bounded domain with smooth boundary. Using the
Lusternik-Schnirelman theory, we relate the number of solutions with
the topology of $\Omega$.
\end{abstract}

\section{Introduction}

In this paper, we are concerned with the multiplicity of nontrivial
solutions for the following class of complex problems

$$
\ \  \left\{
             \begin{array}{l}
             (-i\nabla - A(\mu x))^{2}u = \mu |u|^{q-2}u + |u|^{2^{*}-2}u \ \mbox{in} \ \Omega
\mbox{}\\
        \\ u \in H^{1}_{0}(\Omega,\mathbb{C}),
             \end{array}
           \right.
           \eqno{(P_{\mu})}
$$
where  $\Omega$ is a bounded domain with smooth boundary in
$\mathbb{R}^{N}$, $N\geq 4$, $\mu$ is a positive parameter, $2 \leq q < 2^{*}=\frac{2N}{N-2}$ and $A: \mathbb{R}^{N} \to \mathbb{R}^N$ is a magnetic field belonging to $C(\mathbb{R}^{N},\mathbb{R}^N)\bigcap L^{\infty}(\mathbb{R}^{N},\mathbb{R}^N)$.

\vspace{0,5 cm}

This class of problem is related with the existence of solitary
waves, namely solutions of the form $\psi(x,t) := e^{-i
\frac{E}{h}t}u(x)$, with $E \in \mathbb{R}$, for the nonlinear
Schr\"odinger equation
$$
i h \dfrac{ \partial\psi}{ \partial t } = \left( \dfrac{h}{i}
\nabla - A(z) \right)^2 \psi + U(z)\psi - f(|\psi|^2)\psi,~~z \in
\Omega, \eqno{(NLS)}
$$
where $t >0$, $N \geq 2$, $h$ is the Planck constant and $A$ is a magnetic potential associated to a given
magnetic $B$, $U(x)$ is a real electric potential and the nonlinear
term $f$ is a superlinear function. A direct computation shows
that $\psi$ is a solitary wave for $(NLS)$ if, and only if, $u$ is
a solution of the following problem
\begin{equation} \label{pe}
\left( \dfrac{h}{i}\nabla - A(z)\right)^2 u + V(z)u =
f(|u|^2)u,~~\mbox{ in } \Omega,
\end{equation}
where $V(z) = U(z)-E$. It is important  to investigate the
existence and the shape of such solutions in the semiclassical
limit, namely, as $h \to 0^+$. The importance of this study relies
on the fact that the transition from Quantum Mechanics to
Classical Mechanics can be formally performed by sending the
Planck constant to zero.

There  is a vast literature concerning the existence and
multiplicity of bound state solutions for $(\ref{pe})$ with no
magnetic field, namely $A \equiv 0$ and $h=1$, which becomes an elliptic equation like
$$
\ \  \left\{
             \begin{array}{l}
             - \Delta{u} = \mu |u|^{q-2}u + |u|^{2^{*}-2}u \ \mbox{in} \ \Omega
\mbox{}\\
        \\ u=0 ~~~~ \mbox{on} ~~~~ \partial \Omega.
             \end{array}
           \right.
           \eqno{(P)}
$$

Problem $(P)$ has received considerable attention in last years, after the seminal paper due to Brezis and Nirenberg \cite{BN},
who investigated $(P)$ in the case $q = 2$.  Motivated by that article, many authors have also considered  a lot
of problems involving critical growth in bounded and unbounded domains, see, for example,
Struwe \cite{Struwe}, Garcia Azorero and Peral Alonso \cite{GP, GP2}, Bahri and Coron \cite{BaCo},
Rey \cite{Rey}, Benci and Cerami \cite{BCerami1,BCerami2,BCerami3,BCerami4}, Coron \cite{Coron}, Alves and Ding
\cite{AD}, Alves \cite{Alves} and references therein. This class of problem aroused the
interest of all due to the lack of compactness in the inclusion of
$$
H^{1}_0(\Omega) \hookrightarrow L^{2^{*}}(\Omega),
$$
hence, the associated energy functionals do not satisfy in general the
Palais-Smale condition.

Multiplicity of solutions to $(P)$ involving the geometry of $\Omega$, precisely, the Lusternik�Schnirelman category
$cat_{\Omega}(\Omega)$, was proved in \cite{Rey} for $N \geq 5$ and in \cite{Lazzo} for $N = 4$, cf. \cite{Willem}. Other results of multiplicity involving subcritical growth and category $cat_{\Omega}(\Omega)$ can be found in \cite{BCerami1,BCerami2,CD}. Here, cat$_X(Y)$ denotes the Ljusternik-Schnirelman category of
$Y$ in $X$, namely the least number of closed and contractible sets in the topological space $X$ which cover the closed set $Y \subset X$.

If we  now consider the magnetic case $A \not\equiv 0$, the first result was obtained by Esteban and Lions \cite{EL}.
They have used the concentration-compactness principle and
minimization arguments to obtain solution for $h>0$ fixed and
dimensions $N=2$ or $N=3$. More recently, Kurata \cite{Kurata} proved
that the problem has a least energy solution for any $h>0$ when a
technical condition relating $V$ and $A$ is assumed. Under this
technical condition, he  proved that the associated functional
satisfies the Palais-Smale compactness condition at any level. We
also would like to cite the papers
\cite{CinSec1,CinSec2,ChaSzu,Tan,CinJeaSec,AlvFigFur1,AlvFigFur2} for other
results related to the problem (\ref{pe}) in the presence of
magnetic field.

In view of the results of Rey \cite{Rey} and Lazzo \cite{Lazzo}, it is
natural to ask if the same kind of result holds for the problem
with magnetic field. The main goal of this paper is to present a
positive answer to this question. So, we relate the number of
solution for $(P_\mu)$ with topology of the set $\Omega$ when
the parameter $\mu$ is small. We prove that, for small values
of $\mu$, the magnetic field does not play any role on the
numbers of solutions of the equation $(P_{\mu})$ and therefore
a result in the same spirit of \cite{Rey} and \cite{Lazzo} holds.

\vspace{0,5 cm}

Our main result is:

\begin{theo}\label{gio1}
Let $2\leq q < 2^{*}$. Then, there exists $\mu^{*}>0$ such that, for
each $\mu \in (0, \mu^{*})$, problem $(P_{\mu})$ has at least
$cat_{\Omega}(\Omega)$ nontrivial solutions.
\end{theo}

In the proof of Theorem \ref{gio1}, we apply variational methods and Ljusternik-Schnirelmann
theory.  We follow some arguments developed in \cite{Rey}, \cite{Lazzo} and \cite{AD}, where the non-magnetic case is handled. Is is worthwhile to mention that, since we deal with different problems, where the function are complex, it is necessary to make a careful analysis in some estimates used in that papers.

The paper is organized as follows. In the next section we present
the variational setting of the problem. In Section 3 we prove some preliminary results, and in the Section 4, we prove our main theorem.


\section{Variational framework and notations}


We shall denote by $H^{1}_{0}(\Omega,\mathbb{C})$ the Hilbert space
obtained by the closure of $C^{\infty}_{0}(\Omega, \mathbb{C})$
under the scalar product
$$
\langle u,v \rangle_{A_\mu} := \mbox{Re} \left( \int_{\Omega}
\nabla_{A_\mu}u \overline{\nabla_{A_\mu} v} \; dx \right),
$$
where $A_\mu(x)=A(\mu x)=(A_1(\mu x),A_2(\mu x),...,A_N(\mu x))$, $\mbox{Re}(w)$ denotes the real part of $w \in \mathbb{C}$,
$\overline{w}$ is its complex conjugated, $\nabla_{A_\mu} u
:= (D_{1}u,D_{2}u, ..., D_{N}u )$ and $D_{j}:=-i\partial_{j}-A_j(\mu x
)$, for $j=1,\ldots,N$. The norm induced by this inner product is
given by
$$
\|u\|_{A_\mu}:= \left(\int_{\Omega} |\nabla_{A_\mu}u|^2 dx \right)^{1/2}.
$$

As proved by Esteban and Lions in \cite[Section II]{EL}, for any $u
\in H^{1}_{0}(\Omega,\mathbb{C})$, there holds the \emph{diamagnetic
inequality}, namely
\begin{equation} \label{diamagnetica}
|\nabla |u|(x)|= \left|\mbox{Re}\left(\nabla u
\displaystyle\frac{\overline{u}}{|u|}\right)\right|=\left|\mbox{Re}\biggl((\nabla
u - iA_\mu u)\displaystyle\frac{\overline{u}}{|u|}\biggl)\right|\leq
|\nabla_{A_\mu} u(x)|.
\end{equation}
Thus, if $u \in H^{1}_{0}(\Omega,\mathbb{C})$, we have that $|u|$ belongs to
the usual Sobolev space $H^1_{0}(\Omega,\mathbb{R})$. Moreover, the
embedding $H^{1}_{0}(\Omega,\mathbb{C}) \hookrightarrow
L^q(\Omega,\mathbb{C})$ is continuous for each $1 \leq q \leq 2^*$
and it is compact for $1 \leq q <2^*$.

From now on, we say that a function $u\in H^{1}_{0}(\Omega,\mathbb{C})$ is a weak
solution of $(P_{\mu})$ if
$$
\mbox{Re}\biggl( \int_{\Omega} \nabla_{A_\mu}u
\overline{\nabla_{A_\mu}v}\; dx-\mu \int_{\Omega} |u|^{q-2}u\overline{v}\;dx-
\int_{\Omega} |u|^{2^{*}-2}u\overline{v}\;dx\biggl)=0,
$$
for all $v \in H^1_{0}(\Omega,\mathbb{C})$.

In this paper, the main tool used to prove Theorem \ref{gio1} is the variational method, where the solutions to
$(P_{\mu})$ are obtained  by looking for critical points of the
functional
$$
I_{\mu}(u) =\ds\frac{1}{2}\ds\int_{\Omega}|\nabla_{A_\mu} u|^2 \ \
dx-\frac{\mu}{q}\ds\int_{\Omega}|u|^{q}\ \ dx-
\frac{1}{2^{*}}\ds\int_{\Omega}|u|^{2^{*}}\ \ dx.
$$

A direct computation shows that $ I_{\mu} \in
C^{1}(H^1_{0}(\Omega,\mathbb{C}))$ with
$$
I'_{\mu}(u)v=\mbox{Re}\biggl( \int_{\Omega} \nabla_{A_\mu}u
\overline{\nabla_{A_\mu}v}\;dx- \int_{\Omega} \mu
|u|^{q-2}u\overline{v}\;dx-\int_{\Omega} |u|^{2^{*}-2}u\overline{v}\;dx \biggl).
$$
Thus the weak solutions of $(P_{\mu})$ are  precisely the critical
points of $I_{\mu}$.

\vspace{0.5 cm}

Hereafter, we denote by $\lambda_1 >0$ the best constant of the compact embedding
$$
H^1_{0}(\Omega,\mathbb{C}) \hookrightarrow L^{2}(\Omega,\mathbb{C})
$$
which is given by
$$
\lambda_{1}=\displaystyle\inf_{u \in
H^{1}_{0}(\Omega,\mathbb{C}) \setminus \{0\}}\displaystyle\frac{\int_{\Omega}|\nabla_{A_\mu}
u|^{2} dx }{\left(\int_{\Omega}| u|^{2} dx \right)^{\frac{1}{2}}}.
$$

Moreover, we denote by $S$ the best Sobolev constant of the embedding
$$
H^1_{0}(\Omega,\mathbb{R}) \hookrightarrow L^{2^{*}}(\Omega,\mathbb{R})
$$
which is given by
$$
S=\inf_{u \in
H^{1}_{0}(\Omega,\mathbb{R})\setminus \{0\}}\frac{ \int_{\Omega}|\nabla
u|^{2} dx }{( \int_{\Omega}| u|^{2^{*}} dx )^{2/2^{*}}}.
$$
It is well known that $S$ is independent of $\Omega$ and it is never
achieved, except when $\Omega=\RR^{N}$. Moreover,
$$
S:= \displaystyle\frac{\displaystyle\int_{\mathbb{R}^{N}}|\nabla
U|^{2}\ \ dx}{\biggl(\displaystyle\int_{\mathbb{R}^{N}}|U|^{2^{*}}\
\ dx\biggl)^{2/2^{*}}},
$$
where $U(x)=\displaystyle\frac{C_{N}}{(|x|^{2}+1)^{(N-2)/2}}$ and
$C_{N}$ is a constant such that
$$
-\Delta U=U^{2^{*}-1} \ \mbox{in} \ \ \mathbb{R}^{N}.
$$
A direct computation implies that for all $\epsilon >0$ and $y \in
\mathbb{R}^{N}$ the function
$$
U_{\epsilon, y}(x)=\epsilon^{\frac{2-N}{2}}U(\frac{x-y}{\epsilon})
$$
verifies the equality below
$$
\displaystyle\int_{\mathbb{R}^{N}}|\nabla U_{\epsilon, y}|^{2}\ \
dx=\displaystyle\int_{\mathbb{R}^{N}}|U_{\epsilon, y}|^{2^{*}}\ \
dx= S^{N/2}.
$$

\vspace{0.5 cm}

\begin{lem} \label{SSA} If
$$
S_{A_\mu}=\displaystyle\inf_{u \in
H^{1}_{0}(\Omega,\mathbb{C})}\displaystyle\frac{\displaystyle \int_{\Omega}|\nabla_{A_\mu}
u|^{2} dx }{(\displaystyle \int_{\Omega}| u|^{2^{*}} dx )^{2/2^{*}}},
$$
we have that $S=S_{A_\mu}$.
\end{lem}

\noindent\textbf{Proof.} First of all, we observe that by diamagnetic inequality,
$$
S \leq S_{A_\mu}.
$$
Now, we will prove that $ S \geq S_{A_\mu}$. To this end, we fix $x_0 \in \Omega$. Thus, there exists $r>0$ such that $B_{r}(x_0) \subset \Omega$. Let $\phi$ be a nonnegative smooth cutoff function, such
that
$$
\phi(x)=1 ~~~\mbox{if} ~~~|x|<r, \phi(x)=0 ~~~ \mbox{if} ~~~ |x|>2r,
$$
$$
u_{\epsilon}(x)= \phi(x-x_0)U_{\epsilon, x_0}(x)
$$
and
$$
v_{\epsilon}=\frac{u_{\epsilon}}{|u_{\epsilon}|_{2^{*}}}.
$$

From \cite{BN},
$$
\|v_{\epsilon}\|^{2}= S + O( \epsilon ^{\frac{N-2}{2}})
$$
and
$$
|v_{\epsilon}|^{q}_{q}\rightarrow 0 \,\,\, \mbox{as} \,\,\, \epsilon \to 0 ~~~~ \forall q \in [2 , 2^{*}),
$$
from where it follows that
$$
|v_{\epsilon}|^{2}_{2}\rightarrow 0 \,\,\, \mbox{as} \,\,\, \epsilon \to 0.
$$
From definition of $S_{A_\mu}$, we derive that
$$
S_{A_\mu} \leq \frac{\int_{\Omega}|\nabla_{A_\mu} (e^{i\tau_{x_0}(x)}v_{\epsilon})|^{2} dx }{\left(\int_{\Omega}|e^{i\tau_{x_0}(x)}v_{\epsilon}|^{2^{*}} dx \right)^{\frac{2}{2^{*}}}}
$$
where $\tau_{x_0}(x) := \sum_{j=1}^N A_j(\mu x_0)x^j$. This way,
$$
S_{A_\mu} \leq \frac{\int_{\Omega}|\nabla v_\epsilon|^{2} dx +\int_{\Omega}(A(\mu x_0)-A(\mu x))|v_{\epsilon}|^{2} dx }{\left(\int_{\Omega}|v_{\epsilon}|^{2^{*}} dx \right)^{\frac{2}{2^{*}}}},
$$
or equivalently,
$$
S_{A_\mu} \leq \|v_\epsilon \|^{2}+\int_{\Omega}(A(\mu x_0)-A(\mu x))|v_{\epsilon}|^{2} dx.
$$
Letting $\epsilon \to 0$ and using the fact that $A \in L^{\infty}(\mathbb{R}^{N})$, the above limits leads to
$$
S_{A_\mu} \leq S,
$$
finishing the proof. \hfill\rule{2mm}{2mm}

\vspace{0,5 cm}

\vspace{0.5 cm}

\section {Preliminary results}

Next, we  will show some lemmas related to the functional $I_{\mu}$. Our first lemma is related to the fact that $I_{\mu}$ verifies the mountain pass geometry. However, we omit its proof because it follows by using well known arguments.
\begin{lem}\label{geometriadamontanha}
The functional $I_{\mu}$ satisfies the following conditions:

(i) There exist $\alpha$, $\rho > 0$ such that:
$$
I_{\mu}(u)\geq \alpha \,\,\, \mbox{with} \,\,\,
\|u\|_{A_\mu }=\rho,
$$
for all $\mu>0$ if $2<q< 2^{*}$ and for all $\mu \in
(0,\lambda_{1})$ if $q=2$.

(ii) There exists $e \in B^{c}_{\rho}(0)\subset
H^{1}_{0}(\Omega,\mathbb{C})$ such that $I_{\mu}(e)<0$.
\end{lem}

\vspace*{0.5cm}

Applying the Mountain Pass Theorem without $(PS)$ condition (see
Willem \cite{Willem}), there exists a $(PS)_{b_{\mu}}$ sequence
$(u_{n})\subset H^{1}_{0}(\Omega,\mathbb{C})$, that is, a sequence
satisfying
$$
I_{\mu}(u_{n})\rightarrow b_{\mu} \,\,\, \mbox{and}\,\,\,
I'_{\mu}(u_{n})\rightarrow 0,
$$
where
$$b_{\mu}=\displaystyle\inf_{\gamma \in \Gamma}\displaystyle\max_{t\in [0,1]}I_{\mu}(\gamma(t))
$$
and
$$\Gamma=\{\gamma \in C([0,1],H^{1}_{0}(\Omega,\mathbb{C})): \gamma(0)=0 \;\; \mbox{and} \,\,\,  I_{\mu}(\gamma(1))<0 \}.$$

By standard arguments, $(u_{n})$ is bounded, and so, there exist a
subsequence of $(u_n)$, still denoted by $(u_{n})$, and $u \in
H^{1}_{0}(\Omega,\mathbb{C})$ such that
$$
u_{n}\rightharpoonup u \,\,\, \mbox{in} \,\,\,
H^{1}_{0}(\Omega,\mathbb{C}) \,\,\, \mbox{and}\,\,\,
u_{n}(x)\rightarrow u(x) \,\,\, \mbox{a.e}\,\,\, \mbox{in}\,\,\,
\Omega.
$$

As in \cite[Proposition 3.11]{Rabinowitz}, it is possible to prove
that $b_{\mu}$ verifies the following equalities
$$
b_{\mu}=\widetilde{b}_{\mu}=\widehat{b}_{\mu},
$$
with
$$
\widetilde{b}_{\mu}=\inf\biggl\{\ds\max_{t\geq 0} I_{\mu}(tu) \ : \ u
\in H^{1}_{0}(\Omega,\mathbb{C})\setminus \{0\}\biggl\}
$$
and
$$
\widehat{b}_{\mu}=\inf\biggl\{I_{\mu}(u) \ : \ u \in
\mathcal{N}_{\mu}\biggl\}
$$
where ${\cal N}_\mu$ denotes the Nehari manifold associated with
$I_{\mu}$ given by
$$
{\cal N}_{\mu}=\Big\{u\in H^{1}_{0}(\Omega,\mathbb{C})\setminus
\{0\}:I'_{\mu}(u)u=0\Big\}.
$$

Next, we will prove that $I_{\mu}$ satisfies the local Palais Smale condition.

\begin{lem} \label{H1}
Let $(u_{n}) \subset H^{1}_{0}(\Omega,\mathbb{C})$ be a sequence
that $I_{\mu}(u_{n}) \rightarrow c <
\displaystyle\frac{1}{N}S^{N/2}$ and $\|I_{\mu}'(u_{n})\|=o_{n}(1)$.
Then $I_{\mu}$  satisfies the $(PS)_{c}$ condition for all $\mu>0 $
if $q>2$ and for all $\mu \in (0,\lambda_{1})$ if $q=2$.
\end{lem}
\noindent\textbf{Proof.} Let $(u_{n})\subset
H^{1}_{0}(\Omega,\mathbb{C})$ be a  sequence satisfying
$$
I_{\mu}(u_{n})\rightarrow c \ \ \mbox{and} \ \
I_{\mu}'(u_{n})\rightarrow 0.
$$
From a direct calculus, we have that $(u_{n})$ is bounded in
$H^{1}_{0}(\Omega,\mathbb{C})$. Hence, by diamagnetic inequality,
$(|u_{n}|)$ is bounded in $H^{1}_{0}(\Omega,\mathbb{R})$. Then, for
some subsequence, there is $u \in H^{1}_{0}(\Omega,\mathbb{C})$ such
that $u_{n}\rightharpoonup u$ in $H^{1}_{0}(\Omega,\mathbb{C})$. We
claim that
\begin{eqnarray}\label{conv6}
\displaystyle\int_{\Omega}|u_{n}|^{2^{*}} \ dx\rightarrow
\displaystyle\int_{\Omega}|u|^{2^{*}} \ dx.
\end{eqnarray}
In order to prove this claim, we suppose that
$$
|\nabla |u_n||^2 \rightharpoonup |\nabla |u||^2 + \sigma~~\text{ and
}~~|u_n|^{2^{*}} \rightharpoonup |u|^{2^{*}} +
\nu~~~~~\text{(weak$^*$-sense of measures).}
$$
Using the concentration compactness-principle due to Lions (cf.
\cite[Lemma 1.2]{Lio2}), we obtain a countable index set
$\Lambda$, sequences $(x_i) \subset \overline{\Omega}$, $(\sigma_i),
(\nu_i) \subset (0,\infty)$, such that
\begin{equation}
\nu  =  \sum_{i \in \Lambda}\nu_{i}\delta_{x_{i}},~~\sigma\geq \sum_{i
\in \Lambda}\sigma_{i}\delta_{x_{i}}~~\text{ and }~~S
\nu_{i}^{2/2^{*}}\leq \sigma_{i},
 \label{lema_infinito_eq11}
\end{equation}
for all $i \in\Lambda$, where $\delta_{x_i}$ is the Dirac mass at
$x_i \in \overline{\Omega}$.

Now, for every $\varrho>0$, we set $\psi_{\varrho}(x) :=
\psi((x-x_i)/\varrho)$ where $\psi \in
C_0^{\infty}(\mathbb{R}^{N},[0,1])$ is such that $\psi \equiv 1$ on
$B_1(0)$, $\psi \equiv 0$ on $\mathbb{R}^{N} \setminus B_2(0)$ and
$|\nabla \psi|_{\infty} \leq 2$. Since $(\psi_{\varrho}u_n)$ is
bounded in $H^{1}_{0}(\Omega,\mathbb{C})$ and $\psi_{\varrho}$ takes values in $\mathbb{R}$, a direct calculation
shows that
$$
I_{\mu}'(u_n)(\psi_{\varrho}u_n) \to 0
$$
and
$$
\overline{\nabla_{A_\mu}(u_{n}\psi_{\varrho})}= i\overline{u_{n}}\nabla
\psi_{\varrho}+\psi_{\varrho}\overline{\nabla_{A_\mu}u_{n}}.
$$
Therefore,
$$
\displaystyle\int_{\Omega} \psi_{\varrho}|\nabla_{A_\mu}u_n|^2 \ dx +
\mbox{Re}\biggl( \int_{\Omega}i\overline{u_{n}} \nabla_{A_\mu}u_{n}
\overline{\nabla
\psi_{\varrho}}\biggl)=\mu\displaystyle\int_{\Omega}|u_n|^{q}\psi_{\varrho}
\ dx+ \displaystyle\int_{\Omega} \psi_{\varrho}|u_n|^{2^{**}} \ dx +
o_n(1).
$$
It is not difficult to prove that
$$
\displaystyle\lim_{\varrho\rightarrow
0}[\displaystyle\lim_{n\rightarrow \infty}\mbox{Re}\biggl(
\int_{\Omega}i\overline{u_{n}} \nabla_{A_\mu}u_{n} \overline{\nabla
\psi_{\varrho}}\;dx\biggl)] = 0.
$$
This way, by diamagnetic inequality
$$
\displaystyle\int_{\Omega} \psi_{\varrho}|\nabla|u_n||^2 \ dx \leq
\mu\displaystyle\int_{\Omega}|u_n|^{q}\psi_{\varrho} \ dx+
\displaystyle\int_{\Omega} \psi_{\varrho}|u_n|^{2^{**}} \ dx +
o_n(1).
$$

Consequently, using the fact that $u_n \to u$ in $L^m(\Omega, \mathbb{R})$ for all $1 \leq m
<2^{*}$ and  $\psi_{\varrho}$ has
compact support, we can let $n\to\infty$ in the last inequality to obtain

$$
\int_{\Omega} \psi_{\varrho}\textrm{d}\sigma \leq \int_{\Omega}
\psi_{\varrho}\textrm{d}\nu.
$$
Letting $\varrho \to 0$, it follows that $\nu_i \geq \sigma_i$.
Then, from (\ref{lema_infinito_eq11})
\begin{eqnarray} \label{lema_finito_eq2}
\nu_i \geq \frac{1}{N}S^{N/2}.
\end{eqnarray}
Next, we will prove that the inequality found in (\ref{lema_finito_eq2}) cannot occur, and
therefore the set $\Lambda$ is empty. Indeed, arguing by
contradiction, let us suppose that
$\nu_i\geq\displaystyle\frac{1}{N}S^{N/2}$ for some $i \in \Lambda$.
Once that
\begin{eqnarray*}
c &=& I_{\mu}(u_n) - \displaystyle\frac{1}{2}
I_{\mu}'(u_n)u_n  + o_n(1),
\end{eqnarray*}
it follows that
\begin{eqnarray*}
c \geq \displaystyle\frac{1}{N} \displaystyle\int_{\Omega}
|u_n|^{2^{*}} \ dx + o_n(1)\geq \displaystyle\frac{1}{N}
\displaystyle\int_{B_{\varrho}(x_{i})}  \psi_{\varrho}|u_n|^{2^{*}}
\ dx + o_n(1).
\end{eqnarray*}
Letting $n\to\infty$,
$$
c \geq \displaystyle\frac{1}{N}  \sum_{i \in \Lambda}
\psi_{\varrho}(x_i) \nu_i = \displaystyle\frac{1}{N}  \sum_{i \in
\Lambda} \nu_i \geq \displaystyle\frac{1}{N} S^{N/2},
$$
which does not make sense. Hence, $\Lambda$ is empty and the limit below holds
$$
\displaystyle\int_{\Omega}|u_{n}|^{2^{*}} \ dx\rightarrow
\displaystyle\int_{\Omega}|u|^{2^{*}}\ dx.
$$
The last limit implies that
$$
0\leq \|u_{n}-u\|^{2}_{A_\mu}= I_{\mu}'(u_{n})u_{n}-I_{\mu}'(u_{n})u+
o_{n}(1)=o_{n}(1),
$$
showing that $u_n \to u$ in $H^{1}_{0}(\Omega, \mathbb{C})$.
\hfill\rule{2mm}{2mm}

\vspace{0.5 cm}

The next lemma is a key point in our arguments.

\begin{lem} \label{H2}
The level $b_{\mu}$ verifies the inequality
$$
0< b_{\mu} < \displaystyle\frac{1}{N}S^{N/2},
$$
for all $\mu>0 $ if $q>2$ and for all $\mu \in (0,\lambda_{1})$ if $q=2$.
\end{lem}
\noindent\textbf{Proof.} \; In the sequel, we fix $x_0 \in \Omega$ and $w_{\epsilon}(x)= \tau_{x_0}(x)v_{\epsilon}(x)$ for all $x \in \Omega$, where $\tau_{x_0}$ and $v_{\epsilon}$ were given in the proof of Lemma \ref{SSA}. Setting $g: \mathbb{R} \to \mathbb{R}$ by
$$
g(t)= I_{\mu}(w_{\epsilon})
$$
we have that
$$
g(t)=\ds\frac{t^{2}}{2}\ds\int_{\Omega}|\nabla
v_{\epsilon}|^2\;dx + \ds\frac{t^{2}}{2}\ds\int_{\Omega}(A(\mu x_0)-A(\mu x)) v_{\epsilon}^{2}\;dx -\frac{\mu
t^{q}}{q}\ds\int_{\Omega}|v_{\epsilon}|^{q}\;dx-
\frac{t^{2^{*}}}{2^{*}} .
$$
Thus, there is $t_{\epsilon} >0$ such that
$$
g(t_\epsilon)=\max_{t \geq 0}g(t).
$$
A direct computation shows that $(t_{\epsilon})$ is bounded for $\epsilon$ small enough.
Fixing
$$
h(t)= \ds\frac{t^{2}}{2}\ds\int_{\Omega}|\nabla
v_{\epsilon}|^2 \;dx -\frac{\mu
t^{q}}{q}\ds\int_{\Omega}|v_{\epsilon}|^{q}\;dx-
\frac{t^{2^{*}}}{2^{*}}
$$
and repeating the same arguments  explored in \cite{BN}, we obtain
\begin{equation} \label{Y1}
\max_{t \geq 0} h(t) < \displaystyle\frac{1}{N}S^{N/2} ~~~~ \mbox{for} ~~~~ \epsilon \approx 0.
\end{equation}
On the other hand, once that $A$ is a continuous function, $(t_{\epsilon})$ is bounded, and $v_{\epsilon} \to 0$ in $L^{2}(\Omega)$, we have that
\begin{equation} \label{Y2}
\ds\frac{t_{\epsilon}^{2}}{2}\ds\int_{\Omega}(A(\mu x_0)-A(\mu x)) |v_{\epsilon}|^{2}\;dx \to 0 \,\,\, \mbox{as} \,\,\, \epsilon \to 0.
\end{equation}
Combining (\ref{Y1}) and (\ref{Y2}),
$$
g(t_\epsilon)=\max_{t \geq 0}g(t) < \displaystyle\frac{1}{N}S^{N/2}
$$
for $\epsilon$ small enough. Now, from the definition of $b_{\mu}$,
$$
b_{\mu} \leq g(t_{\epsilon}) \,\,\ \forall \epsilon >0,
$$
from where it follows that
$$
0< b_{\mu} < \displaystyle\frac{1}{N}S^{N/2}.
$$
\hfill\rule{2mm}{2mm}

An immediate consequence of Lemmas \ref{H1} and \ref{H2} is the following result.

\begin{theo} \label{H3} On the hypotheses of Lemma \ref{H2}, the mountain pass level $b_{\mu}$ is a critical value of $I_{\mu}$, that is,
there is $u_{\mu} \in H^{1}_{0}(\Omega,\mathbb{C})$ such that
$$
I_{\mu}(u_{\mu})=b_{\mu} \;\;\; \mbox{and} \;\;\; I'(u_{\mu})=0.
$$
\end{theo}

From now on, we denote by $c_{0}$, $c_{\mu}$ and $\mathcal{M}_{0}$,
$\mathcal{M}_{\mu}$ the mountain pass levels and the Nehari
manifolds associated with the functionals
$$
J_{0}(u) =\ds\frac{1}{2}\ds\int_{\Omega}|\nabla u|^2\;dx -
\frac{1}{2^{*}}\ds\int_{\Omega}|u|^{2^{*}}\;dx \,\, \forall u \in
H^{1}_{0}(\Omega,\mathbb{R})
$$
and
$$
J_{\mu}(u) =\ds\frac{1}{2}\ds\int_{\Omega}|\nabla u|^2\;dx -
\ds\frac{\mu}{q}\ds\int_{\Omega} |u|^{q}\;dx -
\frac{1}{2^{*}}\ds\int_{\Omega}|u|^{2^{*}}\;dx \,\, \forall u \in
H^{1}_{0}(\Omega,\mathbb{R}),
$$
respectively.

\begin{lem}\label{c0igual}
The minimax level $c_0$ is equal to $\frac{1}{N}S^{N/2}$.
\end{lem}
\noindent\textbf{Proof.}  See proof in \cite{AD}. {\fim}

\begin{lem}\label{convergesc0}
If $\displaystyle\lim_{n\rightarrow +\infty}\mu_{n}=0$, then
$\displaystyle\lim_{n\rightarrow
+\infty}b_{\mu_{n}}=c_{0}=\frac{1}{N}S^{N/2}$.
\end{lem}
\noindent\textbf{Proof.} By Theorem \ref{H3}, there is $(u_{n})
\subset H^{1}_{0}(\Omega, \mathbb{C})$ such that
$$
I_{\mu_{n}}(u_{n})= b_{\mu_{n}} \;\;\; \mbox{and} \;\;\;
I'_{\mu_{n}}(u_{n})= 0.
$$
Choosing $t_{n}>0$ such that $t_{n}|u_{n}| \in \mathcal{M}_{0}$, we
derive from diamagnetic inequality that
$$
c_{0} \leq J_0(t_{n}| u_{n}|) \leq I_{\mu_{n}}(t_{n} u_{n})+
\frac{\mu_{n}t_{n}^{q}}{q}|u_{n}|^{q}_{q},
$$
and so,
\begin{equation}\label{ok}
c_{0} \leq b_{\mu_{n}} + \frac{\mu_{n}t_{n}^{q}}{q}|u_{n}|^{q}_{q}.
\end{equation}

From Lemma \ref{H2}, we have that $ b_{\mu_{n}} \leq \frac{1}{N}S^{N/2}$ for $n$ large enough. Conquently, a direct computation implies that $(u_n)$ and $(t_n)$ are  bounded sequences. This way, (\ref{ok})
leads to
\begin{equation} \label{H4}
c_{0}\leq \liminf_{n \to \infty} b_{\mu_{n}}.
\end{equation}

Now, from Lemmas  \ref{H2} and \ref{c0igual}, for $n$
sufficiently large
\begin{equation} \label{Cla1}
b_{\mu_{n}} < \frac{1}{N}S^{N/2}=c_{0},
\end{equation}
leading to
\begin{equation} \label{H5}
\limsup_{n \to \infty} b_{\mu_{n}} \leq c_{0}.
\end{equation}
From this, the lemma follows combining (\ref{H4}) with (\ref{H5}). {\fim}


\section{Technical lemmas}


In this section, we recall  some lemmas which are crucial in the
proof of the main theorem. The  next two lemmas are due to Lions
\cite{Lio2}  and can be found in Willem \cite[Lemma 1.40]{Willem}.

\begin{lem}\label{LL} Let $(u_n) \subset D^{1,2}(\mathbb{R}^{N})$ a sequence such that
$$
u_{n}\rightharpoonup v \ \ \mbox{in} \ \ D^{1,2}(\mathbb{R}^{N}),
$$
\begin{equation} \label{Z01}
|\nabla (u_{n}- u)|^{2}\rightharpoonup \sigma\ \ \mbox{in} \ \
\mathcal{M}(\mathbb{R}^{N}),
\end{equation}
\begin{equation} \label{Z02}
|u_{n}- u|^{2^{**}}\rightharpoonup \nu  \ \ \mbox{in} \ \ \mathcal{M}(\mathbb{R}^{N})
\end{equation}
and
$$
u_{n}\rightharpoonup u \ \ \mbox{in} \ \ \mathbb{R}^{N}.
$$
Then,
\begin{equation} \label{Z1}
\| \nu \|^{\frac{2}{2^{*}}} \leq S^{-1}\| \sigma \|,
\end{equation}
\begin{equation} \label{Z2}
\limsup_{n \to +\infty} |\nabla u_{n}|_{2}^{2}= |\nabla
u_{n}|_{2}^{2} + \| \sigma \| + \sigma_{\infty},
\end{equation}
\begin{equation} \label{Z3}
\limsup_{n \to +\infty} | u_{n}|_{2^{*}}^{2^{*}}= |
u_{n}|_{2^{*}}^{2^{*}} + \| \nu \| + \nu_{\infty},
\end{equation}
and
\begin{equation} \label{Z4}
\nu_{\infty}^{\frac{2}{2^{*}}} \leq S^{-1}\sigma_{\infty},
\end{equation}
with
$$
\sigma_{\infty}= \lim_{R \to \infty}\limsup_{n \to +\infty}\int_{|x|\geq R}|\nabla u_n|^2\;dx
$$
and
$$
\nu_{\infty}= \lim_{R \to \infty}\limsup_{n \to
+\infty}\int_{|x|\geq R}|u_n|^{2^{*}}\;dx.
$$
Moreover, if $u=0$ and $\| \nu \|^{\frac{2}{2^{*}}}=S^{-1}\| \sigma
\|$, the measures $\nu$ and $\sigma$ are concentrated at a single
point.
\end{lem}

\begin{lem}\label{lions}
Let $(u_{n})\subset H_{0}^{1}(\Omega,\mathbb{R})$ be a sequence with
$|u_{n}|_{2^{*}}=1$ and $\|u_{n}\|^{2}=S+o_{n}(1)$. Then there
exists a sequence $(y_{n},\lambda_{n}) \subset \mathbb{R}^{N} \times
\mathbb{R}$ such that $v_{n}(x):=\lambda_{n}^{\frac{N-2}{2}}
u_{n}(\lambda_{n}x+y_{n})$ contains a convergent subsequence, still
denoted by itself, such that $v_{n}\rightarrow v \in
D^{1,2}(\mathbb{R}^{N})$, $\lambda_{n}\rightarrow 0$ and
$y_{n}\rightarrow y \in\overline{\Omega}$.
\end{lem}

An immediate consequence of the last lemma is the following corollary

\begin{coro} \label{C1} Let $(u_{n})\subset H_{0}^{1}(\Omega,\mathbb{R})$ a sequence with
$$
u_n \in {\cal M}_0 ~~~~ \mbox{and} ~~~~ J_0(u_n) \to c_0.
$$
Then there exists a sequence $(y_{n},\lambda_{n}) \subset \mathbb{R}^{N} \times
\mathbb{R}$ such that $v_{n}(x):=\lambda_{n}^{\frac{N-2}{2}}
u_{n}(\lambda_{n}x+y_{n})$ contains a convergent subsequence, still
denoted by itself, such that $v_{n}\rightarrow v \in
D^{1,2}(\mathbb{R}^{N})$, $\lambda_{n}\rightarrow 0$ and
$y_{n}\rightarrow y \in\overline{\Omega}$.
\end{coro}

\vspace{0.5 cm}

Since $\Omega$ is a smooth bounded domain, we choose $r >0$ small enough so that
$$
\Omega_{r}^{+}=\{x\in \mathbb{R}^{N}: dist (x, \Omega) <r\}
$$
and
$$
\Omega_{r}^{-}=\{x\in \mathbb{R}^{N}: dist (x, \partial\Omega) >r\}
$$
are homotopically equivalent to $\Omega$.

\vspace{0.5 cm}

From now on, we consider the functional $J_{\mu, B_{r}}:
H^{1}_{rad}(B_{r}(0),\mathbb{R}) \to \mathbb{R}$ given by
$$
J_{\mu, B_{r}}(u) =\ds\frac{1}{2}\ds\int_{B_{r}(0)}|\nabla u|^2\;dx
-\frac{\mu}{q}\ds\int_{B_{r}(0)}|u|^{q}\;dx-
\frac{1}{2^{*}}\ds\int_{B_{r}(0)}|u|^{2^{*}}\;dx
$$
where
$$
H^{1}_{rad}(B_{r}(0),\mathbb{R})=\{u \in H^{1}_{0}(B_r(0),\mathbb{R}) : u
\ \mbox{is radial}\}.
$$
Moreover, we denote by $m(\mu)$ the mountain pass level associated
with $J_{\mu, B_{r}}$, which can be characterized by
$$
m(\mu):= \inf \{J_{\mu,B_{r}}(u): u \in \mathcal{M}_{\mu ,B_{r}}\}
$$
where
$$
{\cal M}_{\mu, B_{r}}=\Big\{u\in
H^{1}_{rad}(B_{r}(0),\mathbb{R})\setminus \{0\}:J'_{\mu,
B_{r}}(u)u=0\Big\}.
$$

It is not difficult  to check that Lemmas \ref{H1} and \ref{H2} also hold
for $J_{\mu, B_{r}}$. This way, Theorem \ref{H3} is also
true for $J_{\mu, B_{r}}$, from where it follows that there is a
radial function $v_{\mu} \in {\cal M}_{\mu, B_{r}}$ satisfying
$$
J_{\mu,B_{r}}(v_{\mu})=m(\mu) \;\;\; \mbox{and} \;\;\; J'_{\mu,
B_{r}}(v_{\mu})=0.
$$

\begin{lem} \label{H6}
The level $m(\mu)$ converges to $c_0=\frac{1}{N}S^{N/2}$ as $\mu \to
0$, that is, \linebreak $\displaystyle \lim_{\mu \to
0}m(\mu)=c_0=\frac{1}{N}S^{N/2}$.
\end{lem}
\noindent\textbf{Proof.}  See proof in \cite{AD}.  {\fim}

\vspace{0.5 cm}

In what follows, we fix the map $\Psi:\Omega^{-}_{r}\rightarrow {\cal N}_{\mu}$ given by
$$
\Psi_{\mu}(y)(x)=
\left\{
\begin{array}{l}
t_{\mu,y}e^{i\tau_{y}(x)}v_{\mu} (|x-y|) \ \ \mbox{if} \ \ x \in B_{r}(0) \\
0 \ \ \mbox{otherwise},
\end{array}
\right.
$$
where $\tau_{y}(x) := \sum_{j=1}^N A_j(\mu y)x^j$ and $t_{\mu,y} \in (0,+\infty)$ is such that
$$
t_{\mu,y}\textrm{e}^{i\tau_y(.)}v_{\mu}(|\cdot -y|) \in
 \mathcal{N}_{\mu}.
$$
Moreover, following the notation used in \cite{AD}, we denote
by \linebreak $\beta:\mathcal{N}_{\mu}\rightarrow \mathbb{R}^{N}$
the barycenter function given by
$$
\beta(u):= \frac{1}{\displaystyle\int_{\Omega}|u|^{2^{*}}\;dx}\displaystyle\int_{\Omega}x|u|^{2^{*}}\;dx.
$$

Since $v_{\mu}$ is radial, for each $y \in \Omega^{-}_{r}$,
\begin{eqnarray*}
(\beta\circ\Psi)(y)&=&\frac{1}{\displaystyle\int_{\Omega}v_{\mu}(|x-y|)^{2^{*}}\;dx}\displaystyle\int_{\Omega}xv_{\mu}(|x-y|)^{2^{*}}
\;dx = y.
\end{eqnarray*}

\begin{lemma} \label{NovoLema}
Uniformly for $y \in \Omega^{-}_{r}$, there holds
$$
\lim_{\mu \to
0}I_{\mu}(\Psi_{\mu}(y))=c_{0}.
$$
\end{lemma}

\noindent\textbf{Proof.}  Given two sequences
$\mu_n \to 0$ and $(y_n) \subset \Omega_{r}^-$,
we shall prove that
$$
I_{\mu_n}(\Psi_{\mu_n}(y_n)) \to c_{0} ~~~~ \mbox{as} ~~~~ n \to +\infty.
$$

Let $t_n := t_{\lambda_n,y_n}$ and $v_n=v_{\mu_{n}}$ be as in the definition of
$\Psi_{\mu}$. Using the diamagnetic inequality, we have
$$
m(\mu_{n}) \leq I_{\mu_n}(\Psi_{\mu_n}(y_n))
$$
On the other hand,
$$
I_{\mu}(\Psi_{\mu_n}(y_n)) \leq m(\mu_{n}) + \dfrac{t_{n}^2}{2}
\dis\int_{B_{r}(y)} |A(\mu_n y_n)-A(\mu_n x)||v_n|^2\;dx
$$
from where it follows that
\begin{equation} \label{ZZ1}
I_{\mu}(\Psi_{\mu_n}(y_n)) \leq m(\mu_n) + C \dfrac{t_{n}^2}{2} \left( \displaystyle \int_{\mathbb{R}^{N}}|A(\mu_n y_n)-A(\mu_n x)|^{\frac{N}{N-2}}|v_n|^{2^{*}}\;dx\right)^{\frac{N-2}{N}}.
\end{equation}
A direct computation implies that $(t_n)$ is bounded, hence
$$
I_{\mu_n}(\Psi_{\mu_n}(y_n)) \leq m(\mu_n) + C_1 \left( \displaystyle \int_{\mathbb{R}^{N}}|A(\mu_n y_n)-A(\mu_n x)|^{\frac{N}{N-2}}|v_n|^{2^{*}}\;dx\right)^{\frac{N-2}{N}}.
$$

From Corollary \ref{C1}, there exist $(\lambda_{n})\subset \mathbb{R}$ and
$(z_{n})\subset \mathbb{R}^{N}$ with $\lambda_{n}\rightarrow 0$ and
$z_{n}\rightarrow z \in\overline{\Omega}$, such that the sequence
$h_{n}(x):=\lambda_{n}^{\frac{N-2}{2}} v_{n}(\lambda_{n}x+z_{n})$
contains a convergent subsequence, still denoted by itself, that is,
$$
h_{n}\rightarrow h \;\;\; \mbox{in} \;\;\; D^{1,2}(\mathbb{R}^{N},\mathbb{R})
$$
for some $h \in D^{1,2}(\mathbb{R}^{N},\mathbb{R})$. Using the above notations,
\begin{equation} \label{ZZ2}
I_{\mu_n}(\Psi_{\mu_n}(y_n)) \leq  m(\mu_n) + C_1 \left( \displaystyle \int_{\mathbb{R}^{N}}|A(\mu_n y_n)-A(\mu_n  \lambda_n x + \mu_n z_n)|^{\frac{N}{N-2}}|h_n|^{2^{*}}\;dx\right)^{\frac{N-2}{N}}
\end{equation}
Once that $A$ is continuous and belongs to $L^{\infty}(\mathbb{R}^{N})$, it follows that
\begin{equation} \label{ZZ3}
 \int_{\mathbb{R}^{N}}|A(\mu_n y_n)-A(\mu_n \lambda_n x + \mu_n z_n)|^{\frac{N}{N-2}}|h_n|^{2^{*}}\;dx  \to 0.
\end{equation}
Combining (\ref{ZZ1}), (\ref{ZZ2}) and  (\ref{ZZ3}) with the limit $m(\mu_n) \to c_0$, we derive that
$$
I_{\mu_n}(\Psi_{\mu_n}(y_n)) \to c_0,
$$
finishing the proof. {\fim}

\vspace{0.5 cm}

Given $y \in \Omega_{\lambda}^-$, we have that $\Psi_{\mu}(y)
\in \mathcal{M}_{\mu}$. Moreover, setting
\begin{equation} \label{def-g}
g(\mu) := |I_{\mu}(\Psi_{\mu}(y)) - c_{0}|,
\end{equation}
we have that $g(\mu) \to 0$ as $\mu \to 0$ and $I_{\mu}(\Psi_{\mu}(y)) - c_{0} \leq
g(\mu)$.
Hence, the set
$$
\mathcal{O}_{\mu}:=\{ u \in \mathcal{M}_{\mu}:
I_{\mu}(u)\leq c_{0} +g(\mu)\}
$$
contains the function $\Psi_{\mu}(y)$, showing that $\mathcal{O}_{\mu} \not = \emptyset$.

\vspace{0.5 cm}

\begin{lem}\label{vaifecho}
There exists $\mu^{*}>0$ such that, if $\mu \in (0, \mu^{*})$ and $u
\in \mathcal{O}_{\mu}$, then $\beta(u)
\in \Omega_{r}^{+}$.
\end{lem}
\noindent\textbf{Proof.} Suppose by contradiction that there exist
$\mu_{n}\rightarrow 0$, $u_{n} \in \mathcal{N}_{\mu_{n}}$ and
$I_{\mu_{n}}(u_{n})\leq c_0 + g(\mu_n)$ such that $\beta(u_{n})$ does
not belong to $\Omega_{r}^{+}$.

From diamagnetic inequality, there is $t_n \in [0,1]$ such that $v_n := t_n |u_n| \in
\mathcal{M}_{0}$. Hence,
$$
c_{0} \leq J_{0}(t_n |u_n| ) \leq I_{\mu_n}(t_n u_n) +
\frac{\mu_{n}t_{n}}{q}\ds\int_{\Omega}|u_{n}|^{q} \ dx \leq
I_{\mu_n}(u_n) + o_{n}(1)\leq c_0 + o_{n}(1),
$$
and so,
$$
v_n \in \mathcal{M}_{0},~~\beta(v_n) = ~~\beta(u_n) \not\in
\Omega_{r}^+ ~~~ \mbox{and} ~~~~ J_{0}(v_n) \rightarrow c_{0}.
$$

Using Corollary \ref{C1}, there exist $(\lambda_{n})\subset \mathbb{R}$ and
$(y_{n})\subset \mathbb{R}^{N}$ with $\lambda_{n}\rightarrow 0$ and
$y_{n}\rightarrow y \in\overline{\Omega}$, such that the sequence
$h_{n}(x):=\lambda_{n}^{\frac{N-2}{2}} v_{n}(\lambda_{n}x+y_{n})$
contains a convergent subsequence, still denoted by itself, that is,
$$
h_{n}\rightarrow h \;\;\; \mbox{in} \;\;\; D^{1,2}(\mathbb{R}^{N},\mathbb{R})
$$
for some $h \in D^{1,2}(\mathbb{R}^{N},\mathbb{R}) \setminus \{0\}$. Fixing $\phi \in
C^{\infty}_{0}(\mathbb{R}^{N},\mathbb{R})$ with $\phi(x)=x$ for all
$x \in \overline{\Omega}$, a simple computation gives
$$
\beta(v_{n})=\frac{\int_{\mathbb{R^{N}}}\phi(x)|v_{n}(x)|^{2^{*}}dx}{\int_{\mathbb{R^{N}}}|v_{n}(x)|^{2^{*}}dx},
$$
or equivalently
\begin{equation} \label{beta1}
\beta(v_{n})=\frac{\int_{\mathbb{R}^{N}}\phi(\lambda_{n}x+y_{n})|h_{n}(x)|^{2^{*}}dx}{{\int_{\mathbb{R^{N}}}|h_{n}(x)|^{2^{*}}dx}}.
\end{equation}
Letting $n \to +\infty$, we get
$$
\beta(v_{n}) \to y \in \overline{\Omega},
$$
which is a contradiction. {\fim}

\vspace{0.5 cm}


\section{Proof of Theorem \ref{gio1}}


By a direct computation, there exists $C>0$ such that
\begin{equation} \label{prop-nehari}
\|u\|_{A_\mu} \geq C \,\,\, \forall u \in \mathcal{N}_{\mu}.
\end{equation}
Since we are intending to consider the functional $I_{\mu}$
constrained to $\mathcal{N}_{\mu}$, we will need of the following
result.

\begin{lem}\label{gio17}
The functional $I_{\mu}$ constrained to ${\cal{N}}_{\mu}$ satisfies
the $(PS)_{c}$ condition with $c < \frac{1}{N}S^{N/2}$ for $\mu >0$
if $q>2$ and $\mu \in (0, \mu^{*})$ for $q=2$.
\end{lem}
\noindent\textbf{Proof.} Let $(u_n)$ be a (PS)-sequence for
$I_{\mu}$ constrained to $\mathcal{N}_{\mu}$. Then
$I_{\mu}(u_{n})\rightarrow c$ and
\begin{eqnarray}\label{contra}
I'_{\mu}(u_{n}) = \theta_{n} G_{\mu}'(u_{n}) + o_{n}(1),
\end{eqnarray}
for some $(\theta_{n}) \subset \mathbb{R}$, where
$G_{\mu}:H^{1}_{0}(\Omega,\mathbb{C}) \rightarrow \mathbb{R}$ is
given by
\begin{eqnarray*}
G_{\mu}(v) := \displaystyle\int_{\Omega}|\nabla_{A_\mu} v|^{2} \,dx -\mu
\displaystyle\int_{\Omega}|v|^{q} \,dx-
\displaystyle\int_{\Omega}|v|^{2^{*}}\,dx.
\end{eqnarray*}
We recall that $G'_{\mu}(u_{n})u_{n}\leq 0$. Moreover, standard
arguments show that $(u_n)$ is bounded. Thus, up to a subsequence,
$G_{\mu}'(u_{n}) u_{n} \rightarrow l\leq 0$. If $l \neq 0$, we infer
from (\ref{contra}) that $\theta_{n}=o_{n}(1)$. In this case, we can
use (\ref{contra}) again to conclude that $(u_{n})$ is a (PS)$_{c}$
sequence for $I_{\mu}$ in $H^{1}_{0}(\Omega,\mathbb{C})$ and
therefore $(u_{n})$ has a strongly convergent subsequence. If
$l=0$, it follows that
$$
\displaystyle\int_{\Omega}|u_{n}|^{2^{*}} \ dx\rightarrow 0.
$$
Consequently, $\|u_{n}\|_{A_\mu}\rightarrow 0$, obtaining this way a
contradiction with (\ref{prop-nehari}), finishing the proof of the
lemma. {\fim}

\vspace{0.5 cm}

As a consequence of the above arguments, we obtain the following
result.

\begin{coro}\label{17}
If $u$ is a critical point of $I_{\mu}$ constrained to
$\mathcal{N}_{\mu}$, then $u$ is a nontrivial critical point of
$I_{\mu}$ on $H^{1}_{0}(\Omega,\mathbb{C})$.
\end{coro}

\begin{lem}\label{4.5}
If  $\mu^{*}$  is given by Lemma \ref{vaifecho}, then for each $\mu
\in (0, \mu^{*})$, there holds
$$
\emph{cat}_{\mathcal{O}_{\mu}}(\mathcal{O}_{\mu})\geq
\emph{cat}_{\Omega}(\Omega).
$$
\end{lem}

\noindent\textbf{Proof.} Suppose that
$$
\mathcal{O}_{\mu}=\Upsilon_{1}\cup ...\cup \Upsilon_{n},
$$
where $\Upsilon_{j}$, $j=1,\ldots,n$, is closed and contractible
in $\mathcal{O}_{\mu}$. This means that there exists $h_{j}\in
C([0,1]\times \Upsilon_{j},\mathcal{O}_{\mu})$ such that
$$
h_{j}(0,u)=u,\,\,\,h_{j}(1,u)=u_j,\,\,\,\mbox{for each } u\in
\Upsilon_{j},
$$
and some $u_j \in \Upsilon_j$ fixed. Consider the sets $B_j := \gamma^{-1}(\Upsilon_j)$, $j=1,\ldots,n$, which are closed
in $\Omega_{r}^-$ and satisfy
$$
\Omega_{r}^-= B_{1} \cup \cdots \cup B_{n}.
$$

We define the deformation $g_{j}:[0,1]\times B_{j}\rightarrow \Omega_{r}^+$ given
by
$$
g_{j}(t,y)=\beta(h_{j}(t, \gamma(y)))
$$
are well defined. A standard calculation show that these maps are
contractions of the sets $B_j$ in $\Omega_{r}^+$. Hence that
$$
\mbox{cat}_{\Omega}(\Omega)=\mbox{cat}_{\Omega_{r}^+}(\Omega_{r}^-)\leq
n,
$$
and the proposition is proved. {\fim}

\vspace{0.5 cm}
We are now ready to prove our main result.

\vspace{0.5 cm}

\noindent \emph{\bf Proof of Theorem \ref{gio1}.} Arguing as in the
proof of Proposition \ref{gio17}, we can check that $I_{\mu}$
satisfies the $(PS)_{c}$ condition on $\mathcal{N}_{\mu}$ for $c \in
(0, \frac{1}{N}S^{N/2})$. Thus, we can apply standard
Lusternik-Schnirelman theory and Lemma \ref{4.5} to obtain
$\emph{cat}_{\mathcal{O}_{\mu}}(\mathcal{O}_{\mu}) \geq
\mbox{cat}_{\Omega}(\Omega)$ critical points of $I_{\mu}$ restricted
to $\mathcal{N}_{\mu}$. As in Corollary \ref{17}, each one
of these critical points is a critical point of the unconstrained
functional $I_{\mu}$, and therefore, a nonzero weak solution of the
problem $(P_{\mu})$. {\fim}



\addcontentsline{toc}{chapter}{Bibliografia}

\begin{thebibliography}{99}





\bibitem{Alves}{\sc C.O. Alves,} \newblock{ \em  Existence and Multiplicity of Solution for a Class of Quasilinear
Equations, Adv. Non. Studies 5 (2005), 73-87.}



\bibitem{AD} {\sc C. O. Alves} and {\sc Y. H. Ding, }
\newblock {\em Multiplicity of positive solutions to a p-Laplacian
equation involving critical nonlinearity},
\newblock {\em J. Math. Anal. Appl. 279 (2003), 508-521}.




\bibitem{AlvFigFur1}{\sc C.O. Alves, G.M. Figueiredo} and {\sc M.F. Furtado}, \newblock{ \em Multiple solutions for a
nonlinear Schr\"odinger equation with magnetic fields, Comm. in Partial Differential Equations, 36 (2011), 1-22}



\bibitem{AlvFigFur2}{\sc C.O. Alves, G.M. Figueiredo} and {\sc M.F. Furtado}, \newblock{ \em On the number of solutions of NLS equations with magnetics fields in expanding domains, J. Differential Equations 251 (2011), 2534-2548.}




\bibitem{GP} {\sc J. Garcia Azorero} and {\sc I. Peral Alonso},
\newblock {\em Multiplicity of solutions for elliptic problems with critical
exponent or with a nonsymmetric term, Trans. Amer. Math. Soc. , vol
323 n. 2(1991), 877-895.}

\bibitem{GP2}{\sc J. Garcia Azorero} and {\sc I. Peral Alonso,} \newblock{\em  Existence and non-uniqueness for the p-Laplacian: Nonlinear eigenvalues, Comm. Partial Differential Equations 12 (1987), 1389-1430.}

\bibitem{BaCo}{\sc A. Bahri}  and {\sc J.M. Coron,}
\emph{On a nonlinear elliptic equation involving the critical Sobolev exponent: The effect of
the topology of the domain, Comm. Pure Appl. Math. 41 (1988), 253-294.}


\bibitem{BCerami1}{\sc V. Benci}, and {\sc G. Cerami, }
\newblock {\em The effect of the domain topology on the number of positive solutions of nonlinear
elliptic problems, Arch. Rational Mech. Anal. 114 (1991) 79-93. }


\bibitem{BCerami2}
{\sc V. Benci} \, and \, {\sc G. Cerami, }
\newblock {\em Positive solutions of some nonlinear elliptic problems in exterior
domains. Arch. Rational Mech. Anal. 99 (1987), 283-300 }


\bibitem{BCerami3}
{\sc V. Benci} \, and \, {\sc G. Cerami, }
\newblock {\em Multiple positive solutions of some elliptic problems via the
Morse theory and the domain topology , Cal. Var. 02 (1994), 29-48. }


\bibitem{BCerami4} {\sc V. Benci} and {\sc G. Cerami,}
\newblock{ \em Existence of positive solutions of the equation $\Delta u+a(x)u=|u|^{2^{*}-2}u$ in $\mathbb{R}^{N}$,
J. Funct. Anal. 88 (1990), 90-117.}





\bibitem{BN}
{\sc H. Brezis} \, and \, {\sc L. Nirenberg}, \newblock{\em Positive
solutions of nonlinear elliptic equations involving critical Sobolev
exponents, Comm. Pure Appl. Math. 36 (1983), 437-477}.


\bibitem{CD}
{\sc G. Cerami} \, and \, {\sc D. Passaseo, }
\newblock {\em Existence and multiplicity of positive solutions for nonlinear elliptic problems in
exterior domains with "Rich" topology, Nonlinear Anal. 18 (1992), 109-119. }

\bibitem{ChaSzu}{\sc J. Chabrowski} and {\sc A. Szulkin,} \newblock{ \em On the Schr\"odinger equation involving a
critical Sobolev exponent and magnetic field, Top. Meth. Nonlinear Anal. 25 (2005), 3-21.}


\bibitem{CinJeaSec}{\sc S. Cingolani, L. Jeanjean}  and {\sc S. Secchi,}\newblock{ \em Multi-peak solutions for magnetic
NLS equations without non-degeneracy conditions, ESAIM Control Optim.
Calc. Var. 15 (2009), 653-675.}


\bibitem{CinLa}{\sc S. Cingolani} and {\sc M. Lazzo,} \newblock{\em Multiple semiclassical standing waves for a
class of nonlinear Schr\"odinger equations, Topol. Methods Nonlinear Anal. 10 (1997), 1-13.}

\bibitem{CinSec1}{\sc S. Cingolani} and {\sc S. Secchi,} \newblock{\em Semiclassical limit for nonlinear Schr\"odinger
equations with electromagnetic fields, J. Math. Anal. Appl. 275 (2002),
108-130.}

\bibitem{CinSec2}{\sc S. Cingolani} and {\sc S. Secchi,} \newblock{\em Semiclassical states for NLS equations with
magnetic potentials having polynomial growths, J. Math. Phys. 46 (2005), no. 5, 053503, 19 pp.}


\bibitem{Coron}
{\sc J. M. Coron, },
\newblock {\em Topologie en cas limite des injection de Sobolev. C. R. Acad. Sci. Paris
S\'er. I Math. 299 (1984), 209-212.}

\bibitem{EL} {\sc M. Esteban} and {\sc P.L. Lions}, \newblock{ \em Stationary solutions of
    nonlinear Schr\"{o}dinger equations with an external magnetic
    field, PDE and Calculus of Variations, Vol. I, 401-449,
Progr. Nonlinear Differential Equations Appl. 1, Birkh\"{a}user
Boston, MA, 1989. }



\bibitem{Lazzo}{\sc M. Lazzo,} \newblock{ \em Solutions positives multiples pour une \'equation elliptique non lin\'eaire avec l'exposant critique de Sobolev, C. R. Acad. Sci. Paris 314 (1992), 61-64.}



\bibitem{Lio2} {\sc P.L. Lions},
\newblock{ \em The concentration-compactness principle in the calculus of
variations. The limit case, Rev. Mat. Iberoamericana \textbf{1}
(1985), 145-201.}



\bibitem{Kurata} {\sc K. Kurata,} \newblock{\em  Existence and semi-classical limit of the least energy solution
to a nonlinear Schr\"odinger equation with electromagnetic fields, Nonlinear
Anal. 41 (2000), 763-778.}


\bibitem{Rabinowitz}
{\sc P. H. Rabinowitz, }
\newblock {\em On a class of nonlinear Schr\"{o}dinger equations},
\newblock {\em Z. Angew Math. Phys. 43 (1992), 27-42}.




\bibitem{Rey} {\sc O. Rey,} \newblock{\em A multiplicity result for a variational problem with lack of compactness, Nonlinear Anal. 13 (1989), 1241-1249.}


\bibitem{Struwe} {\sc M. Struwe,} \newblock{\em  A global compactness result for elliptic boundary value problems involving limiting nonlinearities,
Math. Z. 187 (1984), 511-517.}

\bibitem{Tan}{\sc Z. Tang,}\newblock{ \em  Multi-bump bound states of nonlinear Schr\"odinger equations with
electromagnetic fields and critical frequency, J. Differential Equations 245 (2008), 2723-2748.}

\bibitem{Willem} {\sc M. Willem,} \newblock {\em Minimax Theorems}, \newblock {\em Birkh\"{a}user, 1996}.


\end{thebibliography}
\end{document}